\documentclass{article}
\usepackage{amssymb}
\usepackage{amsmath}
\usepackage[normalem]{ulem}
\usepackage[utf8]{inputenc}
\usepackage{graphicx}
\usepackage{color}

\newtheorem{theorem}{Theorem}[section]
\newtheorem{lemma}[theorem]{Lemma}
\newtheorem{proposition}[theorem]{Proposition}
\newtheorem{corollary}[theorem]{Corollary}

\newtheorem{preexample}{Example}[section]
\newenvironment{example}{\begin{preexample}}{\end{preexample}}
\newtheorem{preremark}{Remark}
\newenvironment{remark}{\begin{preremark}\rm}{\end{preremark}}

\newenvironment{proof}
{{\bf Proof:}}
{\qquad \hspace*{\fill} $\Box$}%

\title{Projective Controllability of Complete Lifted Control Systems}

\author{ S. N. Stelmastchuk\\
Universidade Federal do Paran\'a, Jandaia do Sul, Brazil}

\begin{document}
	
\maketitle
	
\begin{abstract}
  We introduce the complete lifted control system associated with a control system on a smooth manifold by replacing each vector field with its complete lift to the tangent bundle. We prove that complete lifted control systems are never controllable on the whole tangent bundle, due to the invariance of the zero section. Motivated by this obstruction and by the invariance of complete lifts under fiberwise dilations, we study the induced control system on the projectivized tangent bundle. We establish the relationship between the controllability properties of the lifted and projectivized systems, showing in particular that projective controllability implies controllability of the original system. Our main result provides a sufficient condition for controllability of the projectivized system in terms of a Lie rank condition modulo the Euler vector field. 

\end{abstract}

\section{Introduction}

Let $M$ be a smooth manifold and let $TM$ denote its tangent bundle. Given a vector field $X\in\mathfrak X(M)$, its complete lift is a classical construction in differential geometry that associates to $X$ a vector field $X^c$ on $TM$. The theory of complete lifts has been extensively developed in differential geometry, particularly in the monograph of Yano and Ishihara \cite{yano} and in the survey of Gudmundsson and Kappos \cite{gudmundsson}. From the perspective of geometric mechanics and control, a modern treatment can be found in the work of Bullo and Lewis \cite{bullo}. A fundamental feature of the complete lift is that it describes the infinitesimal evolution of tangent vectors under the differential of a flow. More precisely, if $\varphi_t$ denotes the flow of $X$, then $X^c$ governs the evolution of variations through the tangent map $d\varphi_t$. Thus, while $X$ determines the dynamics on the base manifold, its complete lift encodes the evolution of infinitesimal perturbations along trajectories.

In geometric control theory, the study of variations of trajectories plays a fundamental role in controllability, accessibility and optimal control. This suggests the following natural question: can the complete lift be used to construct a control system on the tangent bundle and, if so, what are its controllability properties?

To address this question, let $X_0,X_1,\ldots,X_m\in\mathfrak X(M)$ and consider the control system
\[
  \dot{x}(t)=X_0(x(t))+\sum_{i=1}^{m}u_i(t)X_i(x(t)),
\]
where $u=(u_1,\ldots,u_m)$ is an admissible control. We define the associated \emph{complete lifted control system} on $TM$ by
\[
  \dot{v}(t)=X_0^c(v(t))+\sum_{i=1}^{m}u_i(t)X_i^c(v(t)).
\]
The resulting dynamics describe simultaneously the evolution of trajectories and the evolution of tangent vectors along them. More precisely, the lifted flow is determined by the differential of the flow of the original system, yielding a coupled dynamics in which the evolution in the fibers is governed by a linear variational equation.

The first objective of this work is to establish the basic properties of complete lifted control systems. In particular, we show that the solutions of the lifted system are naturally expressed in terms of the differential of the flow of the original control system. This characterization reveals an intrinsic geometric obstruction to controllability: the zero section of the tangent bundle is invariant under the lifted dynamics. As a consequence, complete lifted control systems are never controllable on the whole tangent bundle.

A second obstruction arises from the fact that complete lifts commute with the Euler vector field and are therefore invariant under fiberwise dilations. Consequently, the lifted dynamics preserve the radial structure of the fibers and depend essentially on the directions of tangent vectors rather than on their magnitude. This observation suggests that the appropriate framework for studying controllability is not the tangent bundle itself, but rather its projectivization.

The invariance of complete lifts under fiberwise dilations allows the lifted dynamics to descend naturally to the projectivized tangent bundle $P(TM)$. The main objective of this paper is to investigate the controllability properties of the induced control system on $P(TM)$. We establish the relation between the controllability of the complete lifted system and that of its projectivization, showing in particular that projective controllability implies controllability of the original system on $M$. Our main result provides a sufficient condition for controllability of the projectivized system in terms of a Lie rank condition modulo the Euler vector field. Roughly speaking, controllability is obtained whenever the complete lifts generate all directions in the tangent bundle up to the radial direction removed by projectivization.

The paper is organized as follows. In Section~2 we recall the basic properties of complete lifts and fix the notation used throughout the paper. In Section~3 we introduce complete lifted control systems, describe their solutions and establish the fundamental obstruction to controllability on the tangent bundle. Section~4 is devoted to the projectivized dynamics, where we construct the induced control system on $P(TM)$ and prove the main controllability results.

\section{Complete lifts on the tangent bundle}

In this section we recall the definition and the main geometric properties of the complete lift of a vector field. Our presentation follows the classical references of Yano and Ishihara \cite{yano}.

Let $M$ be a smooth manifold of dimension $n$ and let $\pi: TM \to M$ be its tangent bundle. Given a smooth vector field $X \in \mathfrak{X}(M)$, its complete lift $X^c \in \mathfrak{X}(TM)$ is defined as the unique vector field satisfying
\[
  X^c(f^c) = (Xf)^c,
\]
for every smooth function $f \in C^\infty(M)$, where $f^c$ denotes the complete lift of $f$.

A fundamental interpretation of the complete lift is the following: it encodes the differential of the flow of the vector field.

Let $\varphi_t$ be the flow of $X$. Then, for every $w \in T_xM$, one has
\begin{equation}\label{eq:derivative_complete_lift}
  X^c(w) = \left.\frac{d}{dt}\right|_{t=0} (d\varphi_t)_x(w).
\end{equation}
In other words, the complete lift describes the infinitesimal evolution of tangent vectors under the differential of the flow. Thus, while $X$ governs the motion on the base manifold, $X^c$ governs the evolution of velocities along the trajectories of $X$.(See for instance \cite{bullo} for this interpretation).

In local coordinates $(x^i)$ on $M$ and induced coordinates $(x^i, v^i)$ on $TM$, a vector field
\[
  X(x) = X^i(x)\frac{\partial}{\partial x^i}
\]
has complete lift given by
\[
  X^c(x,v) = X^i(x)\frac{\partial}{\partial x^i} + v^j \frac{\partial X^i}{\partial x^j}(x)\frac{\partial}{\partial v^i}.
\]
This expression can be written in a more compact form as
\[
  X^c(x,v) = X(x)\cdot \nabla_x + \big(JX(x)\,v\big)\cdot \nabla_v,
\]
where $JX(x)$ denotes the Jacobian of $X$ at $x$.

This representation reveals a fundamental structural property of the complete lift: its vertical component is linear in the fiber variable $v$. In particular, the evolution of velocities is governed by a linear operator depending on the base point.

Geometrically, this means that the lifted dynamics do not control the magnitude of tangent vectors independently, but only their linear transformation along the trajectory. This feature will be a key ingredient in the analysis of controllability.

Another important property is that the complete lift is $\pi$-related to the original vector field, that is,
\[
  d\pi \circ X^c = X \circ \pi.
\]
Thus, the flow of $X^c$ projects onto the flow of $X$. Moreover, the complete lift preserves the Lie bracket structure:
\[
  [X^c, Y^c] = [X, Y]^c,
\]
for all $X,Y \in \mathfrak{X}(M)$. Consequently, the Lie algebra generated by complete lifts consists again of complete lifts.

A key feature of the complete lift is its invariance under scaling in the fibers. Consider the Euler (or radial) vector field on $TM$ defined by
\[
  R = \sum_{i=1}^n v^i \frac{\partial}{\partial v^i}.
\]
Its flow is given by $(x,v) \mapsto (x, e^s v)$ which represents dilations along the fibers. A direct computation shows that
\[
  [R, X^c] = 0.
\]
Therefore, the complete lift is invariant under the natural scaling action on the tangent bundle.

This invariance has an important geometric consequence: the dynamics induced by $X^c$ depend only on the direction of the tangent vector, and scale linearly with its magnitude. In particular, the radial direction generated by $R$ is not produced by complete lifts.

This observation will be the key obstruction to controllability on $TM$ and will motivate the introduction of the projectivized tangent bundle in Section~4, where the radial direction is factored out.

\begin{example}
Let $X(x) = Ax$ be a linear vector field on $\mathbb{R}^n$, where $A \in \mathbb{R}^{n\times n}$ is a constant matrix. Then $JX(x) = A$, and the complete lift is given by
\[
  X^c(x,v) = (Ax, Av).
\]
Thus, both the position and the velocity evolve linearly under the same operator.
\end{example}

\begin{example}
Let $S^2 \subset \mathbb{R}^3$ be the unit sphere and consider the vector field
\[
  X(x) = Ax,
\]
where $A \in \mathfrak{so}(3)$ is skew-symmetric. Then $X$ generates rotations on $S^2$, and its complete lift is given by
\[
  X^c(x,v) = (Ax, Av).
\]
This means that both the base point and the velocity vector are rotated by the same infinitesimal rotation. This example anticipates the role of rotations in generating directions in the tangent bundle, which will be relevant in the controllability analysis of the projectivized system.
\end{example}

%

\section{Complete lifted control systems}

Let $M$ be a smooth manifold and let $X_0, X_1, \dots, X_m \in \mathfrak{X}(M)$. Consider the control system on $M$ given by
\begin{equation*}\label{controlsystem}
  \dot{x}(t) = X_0(x(t)) + \sum_{i=1}^m u_i(t)\,X_i(x(t)),
\end{equation*}
where $u = (u_1,\dots,u_m)$ is an admissible control. The {\em complete lifted control system} on the tangent bundle $TM$ is defined by
\begin{equation}\label{completecontrolsystem}
  \dot{v}(t) = X_0^c(v(t)) + \sum_{i=1}^m u_i(t)\,X_i^c(v(t)).
\end{equation}
This system describes not only the evolution of the state $x(t)$, but also the evolution of tangent vectors along the trajectory. The next result provides an explicit description of its solutions.

\begin{theorem}\label{thm:solution_complete_control}
  Let $\phi_{t,u}$ be the flow associated with the control system \eqref{controlsystem}. Then, for every fixed admissible control $u$, the curve
  \begin{equation}\label{eq:Gamma}
    \Gamma(t) := \big(\phi_{t,u}(x_0),\, d\phi_{t,u}(x_0)\,v_0\big) \in TM
  \end{equation}
  is an integral curve of the lifted control system \eqref{completecontrolsystem}, that is,
  \[
    \dot{\Gamma}(t)
    = \Big(X_0^c + \sum_{a=1}^m u^a(t)\,X_a^c\Big)\big(\Gamma(t)\big), 
    \qquad \Gamma(0)=(x_0,v_0).
  \]
\end{theorem}
\begin{proof}
  It is a direct application of Proposition B.1 in \cite{barbero}.
\end{proof}

Theorem \ref{thm:solution_complete_control} gives a precise geometric interpretation of the lifted dynamics. While the base trajectory evolves according to \eqref{controlsystem}, the tangent component evolves according to the differential of the flow. In particular, the lifted flow satisfies
\begin{equation}\label{eq:flow_lifted}
  \phi^c_{t,u}(x,v) = \big(\phi_{t,u}(x),\, d\phi_{t,u}(x)\,v\big).
\end{equation}
Thus, the lifted system governs the evolution of velocities along controlled trajectories.

Using the local expression of the complete lift, one obtains the explicit form of the dynamics. Writing $v(t) = (x(t), \xi(t))$, the lifted system is equivalent to
\begin{eqnarray*}
  \dot{x}(t) & = & X_0(x(t)) + \sum_{i=1}^m u_i(t)\,X_i(x(t)),\\
  \dot{\xi}(t) & = & \Big(JX_0(x(t)) + \sum_{i=1}^m u_i(t)\,JX_i(x(t))\Big)\,\xi(t).
\end{eqnarray*}
This shows that, once the base trajectory is fixed, the evolution in the fiber is governed by a linear time-dependent system.

We now recall the notions of orbit and reachable set. For $x \in M$, define
\[
  \mathcal{O}(x) = \{ y \in M : \exists\, t \in \mathbb{R},\, u \in \mathcal{U}, \ \phi_{t,u}(x)=y \},
\]
\[
  \mathcal{R}(x) = \{ y \in M : \exists\, t \ge 0,\, u \in \mathcal{U}, \ \phi_{t,u}(x)=y \}.
\]
Analogously, denote by $\mathcal{O}^c(v)$ and $\mathcal{R}^c(v)$ the orbit and reachable set of the lifted system.

A fundamental obstruction to controllability arises from the invariance of the zero section.

\begin{theorem}\label{thm:not_controllable}
  Let $\dim M > 0$. Then the complete lifted control system is not controllable on $TM$.
\end{theorem}
\begin{proof}
  Let $Z \subset TM$ be the zero section. From \eqref{eq:flow_lifted}, we have
  \[
    \phi^c_{t,u}(x,0) = (\phi_{t,u}(x), 0),
  \]
  which shows that $Z$ is invariant. Therefore, no trajectory starting at $(x,0)$ can reach a point $(y,w)$ with $w \neq 0$. Hence the system is not controllable on $TM$.
\end{proof}

This lack of controllability is intrinsic to the geometry of the complete lift: the system preserves the linear structure of the fibers and does not generate directions transversal to the zero section. This observation will motivate the introduction of projectivized dynamics in the next section.

\section{Projectivization of complete lifted control systems}

The previous section shows that complete lifted control systems are not controllable on the whole tangent bundle $TM$, due to the invariance of the zero section.

Even when restricted to $TM \setminus Z$, an additional obstruction remains. Indeed, as shown in Section~2, complete lifts commute with the Euler (radial) vector field, and therefore preserve the scaling of tangent vectors. Consequently, the dynamics are invariant under scaling and therefore depend only on the direction of tangent vectors, rather than on their magnitude.

This observation suggests that the natural space to study the dynamics of complete lifted systems is not $TM \setminus Z$, but rather the projectivized tangent bundle
\[
  P(TM) = (TM \setminus Z)/\mathbb{R}^*,
\]
where the radial direction is factored out.

In this section, we show that complete lifted control systems induce well-defined control systems on $P(TM)$, and we investigate their controllability properties.

The projectability of complete lifts follows from their invariance under fiberwise scaling.

\begin{proposition}
  Let \(X \in \mathfrak{X}(M)\) and let \(X^c\) denote its complete lift to \(TM\). Let $q: TM \setminus Z \to \mathbb{P}(TM)$ be the canonical projection.  Then \(X^c\) is projectable with respect to \(q\), that is, there exists a unique vector field \(\widetilde{X}\in \mathfrak{X}(\mathbb{P}(TM))\) such that
  \[
    dq_v\big(X^c(v)\big) = \widetilde{X}(q(v))  \quad \text{for all } v\in TM\setminus Z.
  \]
\end{proposition}
\begin{proof}
  Consider the action of \(\mathbb{R}^*\) on \(TM\setminus Z\) given by 
  \[
    \rho_\lambda(v) = \lambda v.
  \]
  It is direct that \(X^c\) is invariant under this action as follows from \eqref{eq:derivative_complete_lift} that 
  \[
    (d\rho_\lambda)_v\big(X^c(v)\big) = X^c(\lambda v), \quad \forall \lambda \neq 0.
  \]
  Hence, \(X^c\) is invariant under the \(\mathbb{R}^*\)-action.  By standard results on quotients by smooth free actions, it follows that \(X^c\) is projectable with respect to \(q\), and therefore induces a vector field \(\widetilde{X}\) on \(\mathbb{P}(TM)\).
\end{proof}

As a consequence, the complete lifted control system induces a well-defined control system on $P(TM)$.

\begin{corollary}
  Let \(X_0,\dots,X_m \in \mathfrak{X}(M)\). Then the complete lifted control system \eqref{completecontrolsystem} is projectable to \(\mathbb{P}(TM)\), defining a control system on \(\mathbb{P}(TM)\).
\end{corollary}

\begin{theorem}\label{thm:projectivized_flow}
  Let \(M\) be a smooth manifold and let \(X_0,\ldots,X_m\in\mathfrak X(M)\). For each fixed admissible control \(u\), let \(\phi_{t,u}\) be the flow of \eqref{controlsystem} on \(M\), and let \(\phi^c_{t,u}\) be the flow of the complete lifted system \eqref{completecontrolsystem} on \(TM\). Then \(\phi^c_{t,u}\) preserves \(TM\setminus Z\) and induces a flow
  \[
    \widehat{\phi}_{t,u}:\mathbb P(TM)\to \mathbb P(TM)
  \]
  satisfying
  \[
    \widehat{\phi}_{t,u}\circ q = q\circ \phi^c_{t,u},
  \]
  where $q$ is the canonical projection. Moreover,
  \[
    \widehat{\phi}_{t,u}(x,[v]) = \left(\phi_{t,u}(x),\left[d\phi_{t,u}(x)v\right]\right).
  \]
\end{theorem}
\begin{proof}
  For each fixed admissible control \(u\), the solution of the complete lifted  system with initial condition \((x,v)\in TM\) is
  \[
    \phi^c_{t,u}(x,v) = \left(\phi_{t,u}(x),d\phi_{t,u}(x)v\right).
  \]
  Since $d\phi_{t,u}(x)$ is a linear isomorphism,  $v\neq 0$ implies $d\phi_{t,u}(x)v\neq 0$. Therefore \(\phi^c_{t,u}\) preserves \(TM\setminus Z\).

  Now let $q$ be the quotient map. Taking \(w=\lambda v\), with \(\lambda\neq 0\), we have
  \[
    d\phi_{t,u}(x)w = \lambda d\phi_{t,u}(x)v \Longrightarrow \left[d\phi_{t,u}(x)w\right]  =\left[d\phi_{t,u}(x)v\right].
  \]
  Thus the rule
  \[
    \widehat{\phi}_{t,u}(x,[v]) :=  \left(\phi_{t,u}(x),\left[d\phi_{t,u}(x)v\right]\right)
  \]
  is well defined. By construction,
  \[
    \widehat{\phi}_{t,u}(q(x,v))  = q(\phi^c_{t,u}(x,v)).
  \]
  That is,
  \[
    \widehat{\phi}_{t,u}\circ q = q\circ \phi^c_{t,u}.
  \]
  Finally, since \(\phi^c_{t,u}\) is the flow of the complete lifted vector field and since each complete lift \(X_i^c\) is projectable with respect to \(q\), the map \(\widehat{\phi}_{t,u}\) is precisely the flow of the induced  projectivized complete lifted control system on \(\mathbb P(TM)\).
\end{proof}

This expression shows that the projectivized dynamics describe the evolution of tangent directions along controlled trajectories. In contrast with the lifted system on $TM$, the radial scaling is eliminated, and only directional information in the tangent bundle is retained.

\begin{proposition}\label{prop:projective_reachable_set}
  Let \(M\) be a smooth manifold and let \(X_0,\ldots,X_m\in\mathfrak X(M)\). Consider the control system \eqref{controlsystem} on \(M\)
  its complete lifted control system on \(TM\), and the induced projectivized system on \(\mathbb P(TM)\). Then, for every \(t>0\) and every  \((x,[v])\in\mathbb P(TM)\),
  \[
    \widehat{\mathcal A}_t(x,[v]) = \left\{ \left(  \phi_{t,u}(x),  \left[d\phi_{t,u}(x)v\right]\right):  u\in\mathcal U\right\}.
  \]
  Equivalently,
  \[
    \widehat{\mathcal A}_t(x,[v]) = q\big(\mathcal A_t^c(x,v)\big).
  \]
\end{proposition}
\begin{proof}
  This result follows directly from the explicit expression of the projectivized flow in above Proposition. 
\end{proof}

Thus, projective controllability allows one to control both the base point and the direction of tangent vectors, but not their magnitude.


\begin{theorem}\label{thm:projective_obstruction}
  Let \(M\) be a smooth manifold and consider the control system \eqref{controlsystem} on \(M\), its complete lifted system on \(TM\setminus Z\), and the induced  projectivized system on \(\mathbb P(TM)\). If the complete lifted system is controllable on \(TM\setminus Z\), then the  projectivized system is controllable on \(\mathbb P(TM)\).
  Conversely, if the projectivized system is controllable on \(\mathbb P(TM)\), then for every \((x,v),(y,w)\in TM\setminus Z\), there exist \(t>0\), an  admissible control \(u\), and \(\lambda\in\mathbb R^*\) such that
  \[
    \phi^c_{t,u}(x,v)=(y,\lambda w).
  \]
  Thus projective controllability controls the base point and the tangent direction, but not necessarily the fiber scale.
\end{theorem}
\begin{proof}
  Let $q$ be the canonical projection. Assume first that the complete lifted system is controllable on  \(TM\setminus Z\). Let \((x,[v]),(y,[w])\in\mathbb P(TM)\). Choose representatives \(v\in T_xM\setminus\{0\}\) and  \(w\in T_yM\setminus\{0\}\). By controllability on \(TM\setminus Z\), there exist \(T>0\) and an admissible control \(u\) such that $\phi^c_{t,u}(x,v)=(y,w)$. Applying \(q\), and using the identity
  \[
    \widehat\phi_{t,u}\circ q=q\circ\phi^c_{t,u},
  \]
  we obtain
  \[
    \widehat\phi_{t,u}(x,[v])=(y,[w]).
  \]
  Hence the projectivized system is controllable on \(\mathbb P(TM)\). 
  
  Conversely, assume that the projectivized system is controllable on  \(\mathbb P(TM)\). Let \((x,v),(y,w)\in TM\setminus Z\). Then there exist \(t>0\) and an admissible control \(u\) such that $\widehat\phi_{t,u}(x,[v])=(y,[w])$. Since
  \[
    \widehat\phi_{t,u}(x,[v]) = \left( \phi_{t,u}(x),\left[d\phi_{t,u}(x)v\right]\right),
  \]
  we get $\phi_{t,u}(x)=y$ and $\left[d\phi_{t,u}(x)v\right]=[w]$. Therefore, there exists \(\lambda\in\mathbb R^*\) such that $d\phi_{t,u}(x)v=\lambda w$. Thus
  \[
    \phi^c_{t,u}(x,v)=(y,\lambda w).
  \]
\end{proof}

Projective controllability removes the radial obstruction. Therefore, the difference between controllability on \(TM\setminus Z\) and controllability on \(\mathbb P(TM)\) lies precisely in the ability, or inability, to control the scaling of tangent vectors.

\begin{remark}
  For complete lifts one has $[X^c,Y^c]=[X,Y]^c$. Consequently, the Lie algebra generated by complete lifts is again made of complete lifts. In particular, it does not generate the Euler radial vector  field $R=\sum_i y^i\frac{\partial}{\partial y^i}$, unless an additional vertical or radial mechanism is added to the system. Thus, for pure complete lifted systems, projective controllability does not generally imply controllability on \(TM\setminus Z\).
\end{remark}

\begin{example}
  Let \(M=\mathbb R^n\) and consider the fully actuated control system
  \[
    \dot x=u,\qquad u\in\mathbb R^n.
  \]
  Then the base system is controllable on \(\mathbb R^n\). However, since the controlled vector fields are constant, their complete lifts are
  \[
    \left(\frac{\partial}{\partial x^i}\right)^c  = \frac{\partial}{\partial x^i}.
  \]
  Hence the complete lifted system is
  \[
    \dot x=u,\qquad \dot y=0.
  \]
  Therefore neither the direction nor the norm of \(y\) can be changed. The system is not controllable on \(T\mathbb R^n\setminus Z\), and its  projectivization is not controllable on \(\mathbb P(T\mathbb R^n)\).  This shows that controllability of the base system alone does not imply projective controllability of the complete lift.
\end{example}

We now show that projective controllability already implies controllability of the original system.

\begin{theorem}\label{thm:projective_implies_base}
  Let \(M\) be a smooth manifold with \(\dim M\geq 1\). Consider the control system \eqref{controlsystem} on \(M\), and let its complete lift induce a projectivized control system on  \(\mathbb P(TM)\). If the projectivized control system is controllable on \(\mathbb P(TM)\), then the original control system is controllable on \(M\).
\end{theorem}

\begin{proof}
  Let \(x,y\in M\). Since \(\dim M\geq 1\), the projective fibers \(\mathbb P(T_xM)\) and \(\mathbb P(T_yM)\) are nonempty. Choose arbitraryclasses
  \[
    [v]\in \mathbb P(T_xM), \qquad  [w]\in \mathbb P(T_yM).
  \]
  By controllability of the projectivized system on \(\mathbb P(TM)\), there exist \(t>0\) and an admissible control \(u\) such that $    \widehat\phi_{t,u}(x,[v])=(y,[w])$. But the projectivized flow is given by
  \[
    \widehat\phi_{t,u}(x,[v]) = \left(\phi_{t,u}(x),\left[d\phi_{t,u}(x)v\right]\right).
  \]
  Therefore,
  \[
    \left(\phi_{t,u}(x),\left[d\phi_{t,u}(x)v\right]\right)=(y,[w]).
  \]
  Comparing the base components, we obtain $\phi_{t,u}(x)=y$. Hence \(y\in\mathcal A(x)\). Since \(x,y\in M\) were arbitrary, the original control system is controllable on \(M\).
\end{proof}


In order to properly understand the structure of the projectivized tangent bundle \(\mathbb{P}(TM)\), it is essential to identify the geometric directions that are eliminated by the projection $q: TM \setminus Z \longrightarrow \mathbb{P}(TM)$. These directions are precisely encoded by the Euler (or radial) vector field on \(TM\). The next result makes this relation precise at the infinitesimal level.

\begin{proposition}\label{prop:euler_kernel}
  Let $q: TM\setminus Z \longrightarrow \mathbb{P}(TM)$ be the canonical projection. Then, for every \(v\in TM\setminus Z\), $\ker dq_v = \operatorname{span}\{R_v\}$.
\end{proposition}
\begin{proof}
  The fibers of \(q\) are given by
  \[
    q^{-1}(q(v)) = \{ \lambda v : \lambda \in \mathbb{R}^* \},
  \]
  which are precisely the orbits of the scaling action. These fibers are one-dimensional submanifolds of \(TM\setminus Z\). Their tangent space at \(v\) is generated by the velocity vector of the curve \(s \mapsto e^s v\), that is,
  \[
    \left.\frac{d}{ds}\right|_{s=0} e^s v = R(v).
  \]
  Therefore,
  \[
    T_v\big(q^{-1}(q(v))\big) = \operatorname{span}\{R_v\}.
  \]
  Since the kernel of the differential of a submersion coincides with the tangent space to the fiber, it follows that $\ker dq_v = \operatorname{span}\{R_v\}$.
\end{proof}


Proposition~\ref{prop:euler_kernel} shows that the Euler vector field spans exactly the directions that are collapsed by the projectivization map. In other words, the radial directions correspond to pure scalings of tangent vectors, and therefore do not carry any directional information in \(\mathbb{P}(TM)\).

This interpretation plays a crucial role in the analysis of control systems on \(\mathbb{P}(TM)\). Indeed, if a family of vector fields on \(TM\setminus Z\) generates all directions except possibly the radial one, then, after projection, it generates the entire tangent space of \(\mathbb{P}(TM)\).

In this sense, the Euler vector field isolates the purely radial dynamics, which are irrelevant for the projectivized system. Once these directions are factored out, the remaining dynamics encode precisely the evolution of directions in the tangent bundle.

We now formulate a sufficient condition for controllability of the projectivized system. The key idea is that it suffices to generate all directions modulo the radial direction spanned by the Euler vector field.

The projective tangent bundle removes exactly one direction at each point, namely the radial direction generated by the Euler field.

\begin{theorem}\label{thm:projective_controllability_euler}
  Let \(M\) be a connected smooth manifold and consider the complete lifted control system on \(TM\setminus Z\), together with its induced system on  \(\mathbb{P}(TM)\). Let
  \[
    R = \sum_{i=1}^n y^i \frac{\partial}{\partial y^i}
  \]
  be the Euler vector field.  Assume that:
  \begin{enumerate}
    \item the projected family of vector fields on \(\mathbb{P}(TM)\) is symmetric;
    \item for every \(v\in TM\setminus Z\),
    \[
      \operatorname{Lie}\{X_0^c,\ldots,X_m^c\}(v) + \operatorname{span}\{R_v\}  = T_v(TM\setminus Z).
    \]
  \end{enumerate}
  Then the projectivized complete lifted control system is controllable on  \(\mathbb{P}(TM)\).
\end{theorem}
\begin{proof}
  Let $q$ be the canonical projection. Since \(q\) is a smooth submersion, for every  \(v\in TM\setminus Z\), the differential
  \[
    dq_v : T_v(TM\setminus Z) \to T_{q(v)}\mathbb{P}(TM)
  \]
  is surjective, and its kernel is given by $\ker dq_v = \operatorname{span}\{R_v\}$. Each complete lifted vector field \(X_i^c\) is projectable with respect to  \(q\), and we denote by \(\widetilde{X}_i\) its projection onto \(\mathbb{P}(TM)\). Since projectability is preserved under Lie brackets, we  have
  \[
    dq_v\Big(\operatorname{Lie}\{X_0^c,\ldots,X_m^c\}(v) \Big)  = \operatorname{Lie}\{\widetilde{X}_0,\ldots,\widetilde{X}_m\}(q(v)).
  \]
  By hypothesis,
  \[
    \operatorname{Lie}\{X_0^c,\ldots,X_m^c\}(v) + \operatorname{span}\{R_v\}  = T_v(TM\setminus Z).
  \]
  Applying \(dq_v\) and using \(dq_v(R_v)=0\), we obtain
  \[
    \operatorname{Lie}\{\widetilde{X}_0,\ldots,\widetilde{X}_m\}(q(v))  = T_{q(v)}\mathbb{P}(TM).
  \]
  Thus, the projected family is bracket-generating on \(\mathbb{P}(TM)\). Since \(M\) is connected, \(\mathbb{P}(TM)\) is also connected. Moreover, the projected family is symmetric, hence attainable sets coincide with  orbits. By Corollary 5.2 in \cite{agrachev}, it follows that the projectivized system is controllable on \(\mathbb{P}(TM)\).
\end{proof}

A natural question is about a condition that caracterize the simmetry on projectable control system. To answer this question we presente the next result. 

\begin{proposition}\label{prop:projective_symmetry_condition}
  Consider the complete lifted control system (\ref{completecontrolsystem}). For \(u\in U\), set $F_u^c = X_0^c+\sum_{i=1}^m u_iX_i^c$. The projected admissible family on \(\mathbb P(TM)\) is symmetric if and only if, for every \(u\in U\), there exists \(v\in U\) such that
  \[
    F_u^c+F_v^c
    \in
    \Gamma(\operatorname{span}\{R\}).
  \]
\end{proposition}
\begin{proof}
  Since \(\pi_*R=0\), two vector fields on \(TM\setminus Z\) induce opposite vector fields on \(\mathbb P(TM)\) precisely when their sum is vertical with respect to the projection \(\pi\), that is, when their sum is tangent to the radial direction generated by \(R\). Therefore,
  $\widehat{F_v^c} = -\widehat{F_u^c}$ if and only if $\pi_*(F_u^c+F_v^c)=0$, which is equivalent to
  \[
    F_u^c+F_v^c
    \in
    \Gamma(\operatorname{span}\{R\}).
  \]
\end{proof}

\begin{corollary}\label{cor:driftless_projective_symmetry}
  If \(X_0=0\) and \(U=-U\), then the projected admissible family on \(\mathbb P(TM)\) is symmetric.
\end{corollary}
\begin{proof}
  In this case,
  \[
    F_u^c=\sum_{i=1}^m u_iX_i^c.
  \]
  Since \(U=-U\), for every \(u\in U\) one has \(-u\in U\), and $F_{-u}^c=-F_u^c$. Hence the projected family is symmetric.
\end{proof}

We now illustrate Theorem~\ref{thm:projective_controllability_euler}  with an example which satisfies the rank condition modulo the  Euler vector field and yields controllability on \(\mathbb{P}(TM)\).
\begin{example}
  Let \(S^2\subset\mathbb R^3\). For \(a\in\mathbb R^3\), denote by
  \[
    Y_a(p)=a\times p,\qquad p\in S^2,
  \]
  the infinitesimal rotation field generated by \(a\). Set $Y_i=Y_{e_i}$, $i=1,2,3$. Consider the control system
  \[
    \dot p
    =
    Y_3(p)+u_1(-Y_3)(p)+u_2Y_1(p)+u_3Y_2(p),
    \quad u=(u_1,u_2,u_3)\in U,
  \]
  where \(U\subset\mathbb R^3\) is symmetric with respect to $\bar u=(1,0,0)$, that is,
  \[
    u\in U
    \quad\Longrightarrow\quad
    2\bar u-u\in U.
  \]
  The drift \(X_0=Y_3\) is nonzero. However,
  \[
    X_0+\bar u_1X_1+\bar u_2X_2+\bar u_3X_3 = Y_3-Y_3=0.
  \]
  Hence $F_{\bar u}^c=0\in\Gamma(\operatorname{span}\{R\})$, and the projected admissible family on \(\mathbb P(TS^2)\) is symmetric.

  We now verify the rank condition. The fields \(Y_1,Y_2,Y_3\) are the infinitesimal generators of the standard action of \(SO(3)\) on \(S^2\).
  Their complete lifts generate the tangent action of \(SO(3)\) on \(TS^2\setminus Z\). For each \(r>0\), this action is transitive on
  \[
    \{(p,v)\in TS^2:\|v\|=r\}.
  \]
  Therefore
  \[
    \dim
    \operatorname{Lie}\{Y_1^c,Y_2^c,Y_3^c\}(p,v)
    =
    3
  \]
  for every \((p,v)\in TS^2\setminus Z\). Since the Euler vector field is
  transverse to these level sets of \(\|v\|\), we obtain
  \[
    \operatorname{Lie}\{X_0^c,X_1^c,X_2^c,X_3^c\}(p,v)
    +
    \operatorname{span}\{R_{(p,v)}\}
    =
    T_{(p,v)}(TS^2\setminus Z).
  \]
  By Theorem~\ref{thm:projective_controllability_euler}, the projectivized
  complete lifted system is controllable on \(\mathbb P(TS^2)\).
\end{example}

This example illustrates that the only obstruction to full controllability in the lifted system is the radial direction. Once this direction is removed, the system becomes fully controllable at the projective level.

\begin{example}
  Let \(G=SO(3)\), and identify its Lie algebra \(\mathfrak{so}(3)\) with \(\mathbb R^3\). Let \(e_1,e_2,e_3\) be the canonical basis of
  \(\mathbb R^3\), and denote by \(X_A\) the right-invariant vector field on \(SO(3)\) determined by \(A\in\mathfrak{so}(3)\), that is, $X_A(g)=dR_gA$. Consider the control system
  \[
    \dot g = X_0(g)+u_1X_1(g)+u_2X_2(g)+u_3X_3(g),
  \]
  where 
  \[
    X_0=X_{e_3},
    \qquad
    X_1=X_{-e_3},
    \qquad
    X_2=X_{e_1},
    \qquad
    X_3=X_{e_2}.
  \]
  Thus the drift \(X_0\) is nonzero.

  Let \(U\subset\mathbb R^3\) be a control set with nonempty interior and symmetric with respect to $\bar u=(1,0,0)$, namely,
  \[
    u\in U \quad\Longrightarrow\quad 2\bar u-u\in U.
  \]
  For \(u\in U\), set
  \[
    F_u=X_0+u_1X_1+u_2X_2+u_3X_3.
  \]
  Since $F_{\bar u}=X_{e_3}+X_{-e_3}=0$, if \(v=2\bar u-u\), then $F_v=-F_u$. Therefore $F_v^c=-F_u^c$, and the projected admissible family on \(\mathbb P(TSO(3))\) is symmetric.

  We now verify the rank condition in Theorem~\ref{thm:projective_controllability_euler}. Under the right trivialization
  \[
    TSO(3)\simeq SO(3)\times\mathfrak{so}(3),
  \]
  the complete lift of \(X_A\) is
  \[
    X_A^c(g,\xi)
    =
    \big(dR_gA,\operatorname{ad}_A\xi\big).
  \]
  Since \(X_{e_1},X_{e_2},X_{e_3}\) belong to the admissible family, the horizontal components of their complete lifts span \(T_gSO(3)\). Their
  vertical components span $[\mathfrak{so}(3),\xi]$. Since the Euler vector field is
  \[
    R_{(g,\xi)}=(0,\xi),
  \]
  it remains to observe that, for every \(\xi\neq0\),
  \[
    [\mathfrak{so}(3),\xi]\oplus\mathbb R\xi
    =
    \mathfrak{so}(3).
  \]
  Indeed, the linear map
  \[
    A\mapsto[A,\xi]
  \]
  has kernel \(\mathbb R\xi\), and hence image of dimension \(2\). Moreover,
  this image has trivial intersection with \(\mathbb R\xi\). Thus
  \[
    [\mathfrak{so}(3),\xi]\oplus\mathbb R\xi
    =
    \mathfrak{so}(3).
  \]
  Consequently,
  \[
    \operatorname{Lie}\{X_0^c,X_1^c,X_2^c,X_3^c\}(g,\xi)
    +
    \operatorname{span}\{R_{(g,\xi)}\}
    =
    T_{(g,\xi)}(TSO(3)\setminus Z)
  \]
  for every \((g,\xi)\in TSO(3)\setminus Z\).
  By Theorem~\ref{thm:projective_controllability_euler}, the projectivized
  complete lifted system is controllable on \(\mathbb P(TSO(3))\).
\end{example}


\begin{thebibliography}{20}

	\bibitem{agrachev} A.~A. Agrachev and Y.~L. Sachkov, {\it Control theory from the geometric viewpoint}, Encyclopaedia of Mathematical Sciences Control Theory and Optimization, 87 II, Springer, Berlin, 2004; MR2062547
  
  \bibitem{barbero} M. Barbero-Linan and M.C. Munoz-Lecanda, {\em Geometric Approach to Pontryaginy Maximum Principle}, Acta Appl Math (2009) 108: 42985.
  
  \bibitem{bullo} F. Bullo and A.~D. Lewis, {\em Geometric Control of Mechanical Systems: Modeling, Analysis, and Design for Simple Mechanical Control Systems}, Springer Science+Business Media, New York, 2005.   
  



  \bibitem{gudmundsson} S. Gudmundsson and E. Kappos, On the geometry of tangent bundles, Expo. Math. {\bf 20} (2002), no.~1, 1--41; MR1888866

  \bibitem{jurdjevic} V. Jurdjevic, {\it Geometric control theory}, Cambridge Studies in Advanced Mathematics, 52, Cambridge Univ. Press, Cambridge, 1997; MR1425878
  
	

  \bibitem{sontag} E.~D. Sontag, {\it Mathematical control theory}, second edition, Texts in Applied Mathematics, 6, Springer, New York, 1998; MR1640001
  
  \bibitem{yano} K. Yano and S. Ishihara, {\it Tangent and cotangent bundles: differential geometry}, Pure and Applied Mathematics, No. 16, Dekker, New York, 1973; MR0350650
  
\end{thebibliography}
\end{document}